\documentclass[reqno]{amsart}

\usepackage{tabu}
\usepackage{amssymb}
\usepackage{mathtools}
\usepackage{a4wide,amsmath}
\usepackage{mathrsfs}
\usepackage{amsthm}
\numberwithin{equation}{section}
\numberwithin{figure}{section}
\numberwithin{table}{section}
\usepackage{bbm}
\usepackage{subfig}
\usepackage{enumerate}
\usepackage[section]{placeins}
\usepackage{graphicx}		  
\usepackage{ifpdf}
\ifpdf
	\DeclareGraphicsExtensions{.pdf,.eps,.jpg,.png}	
	\usepackage[suffix=]{epstopdf}
\fi
\usepackage{xcolor}
\usepackage[utf8]{inputenc}
\usepackage{hyperref}
\hypersetup{hidelinks}

\long\def\MSC#1\EndMSC{\def\arg{#1}\ifx\arg\empty\relax\else
     {\narrower\noindent%
{2020 Mathematics Subject Classification}: #1\\} \fi}
\long\def\PACS#1\EndPACS{\def\arg{#1}\ifx\arg\empty\relax\else
     {\narrower\noindent%
{PACS numbers}: #1}\fi}
\long\def\KEY#1\EndKEY{\def\arg{#1}\ifx\arg\empty\relax\else
	{\narrower\noindent%
Keywords: #1\\}\fi}

\theoremstyle{plain}
\newtheorem{theorem}{Theorem}[section]

\theoremstyle{definition}
\newtheorem{definition}[theorem]{Definition}

\newtheorem{assumption}[theorem]{Assumption}
\theoremstyle{remark}
\newtheorem{remark}[theorem]{Remark}

\newcommand{\inner}[1]{\langle#1\rangle}

\newcommand{\dist}{\mathop{\textup{dist}}}

\newcommand{\rum}[1]{\mathbb{#1}}

\begin{document}

\title{Simplified reconstruction of layered materials in EIT}

\author[H.~Garde]{Henrik Garde}
\address[H.~Garde]{Department of Mathematics, Aarhus University, Ny Munkegade 118, 8000 Aarhus C, Denmark.}
\email{garde@math.au.dk}

\begin{abstract}
	This short note considerably simplifies a reconstruction method by the author (Comm. PDE, 45(9):1118--1133, 2020), for reconstructing piecewise constant layered conductivities (PCLC) from partial boundary measurements in electrical impedance tomography. Theory from monotonicity-based reconstruction of extreme inclusions eliminates most of the bookkeeping related to multiple components of each layer, and also simplifies the involved test operators. Moreover, the method no longer requires a priori lower and upper bounds to the unknown conductivity values.
\end{abstract}

\maketitle

\KEY
electrical impedance tomography, 
partial data reconstruction,
piecewise constant coefficient,
monotonicity principle.
\EndKEY

\MSC
35R30, 
35Q60, 
35R05, 
47H05. 
\EndMSC

\section{Introduction}

This note simplifies the reconstruction method in \cite{Garde_2019b} for piecewise constant layered conductivities, using  the main result from \cite{Garde2020c} on monotonicity-based reconstruction of inclusions. See~\cite{Garde2020c,Garde_2019b} and the references therein for more in-depth information on the methods and on reconstruction in electrical impedance tomography (EIT) in general. The main theoretical developments of monotonicity-based reconstruction in EIT can be found in \cite{Garde2020c,Esposito2021,Garde_2019b,GardeHyvonen2021b,Harrach_2019,Harrach10,Harrach13,Ikehata1998a,Kang1997b,Tamburrino2002} for the continuum model and in \cite{GardeStaboulis_2016,Garde_2019,Harrach_2019,Harrach15} for electrode models. 

Let $\Omega\subset \mathbb{R}^d$, $d\geq 2$, be a bounded Lipschitz domain with connected complement. Let $\nu$ be an outer unit normal to $\partial\Omega$, and let $\Gamma\subseteq\partial\Omega$ be a relatively open subset. 

We consider the partial data conductivity problem, formally written
\begin{equation} \label{eq:condeq}
	-\nabla\cdot(\sigma\nabla u) = 0 \quad\text{in } \Omega, \qquad \nu\cdot(\sigma\nabla u)|_{\partial\Omega} = \begin{cases}
		f &\quad\text{on } \Gamma, \\
		0 &\quad\text{on } \partial\Omega\setminus\Gamma,
	\end{cases}
\end{equation}
where $f$ belongs to
\begin{equation*}
	L^2_\diamond(\Gamma) = \{ g\in L^2(\Gamma) \mid \inner{g,1} = 0\}.
\end{equation*}
Here and in the sequel, $\inner{\cdot,\cdot}$ refers to the usual inner product on $L^2(\Gamma)$. For the electric potential~$u$ we also enforce a $\Gamma$-mean free condition.

A conductivity coefficient $\sigma$ can in general be nonnegative and measurable. We assume that $\sigma$ can formally equal zero or infinity on certain separated Lipschitz sets; such sets are called extreme inclusions. Away from the extreme inclusions, $\sigma$ is assumed to be bounded away from zero and infinity. This gives rise to a local Neumann-to-Dirichlet (ND) map on $\Gamma$, which is a compact self-adjoint operator $\Lambda(\sigma): f\mapsto u|_\Gamma$ on $L^2_\diamond(\Gamma)$. See \cite[Section~2]{Garde2020c} for the precise definitions and assumptions. 

In the following we use the notation that $\gamma$ denotes an unknown conductivity coefficient that we seek to reconstruct from $\Lambda(\gamma)$. We will need several other coefficients in the reconstruction method, and these may take the place of $\sigma$ in the PDE problem~\eqref{eq:condeq} in order to define their local ND maps.

\begin{figure}[htb]
	\centering
	\includegraphics[width=\textwidth]{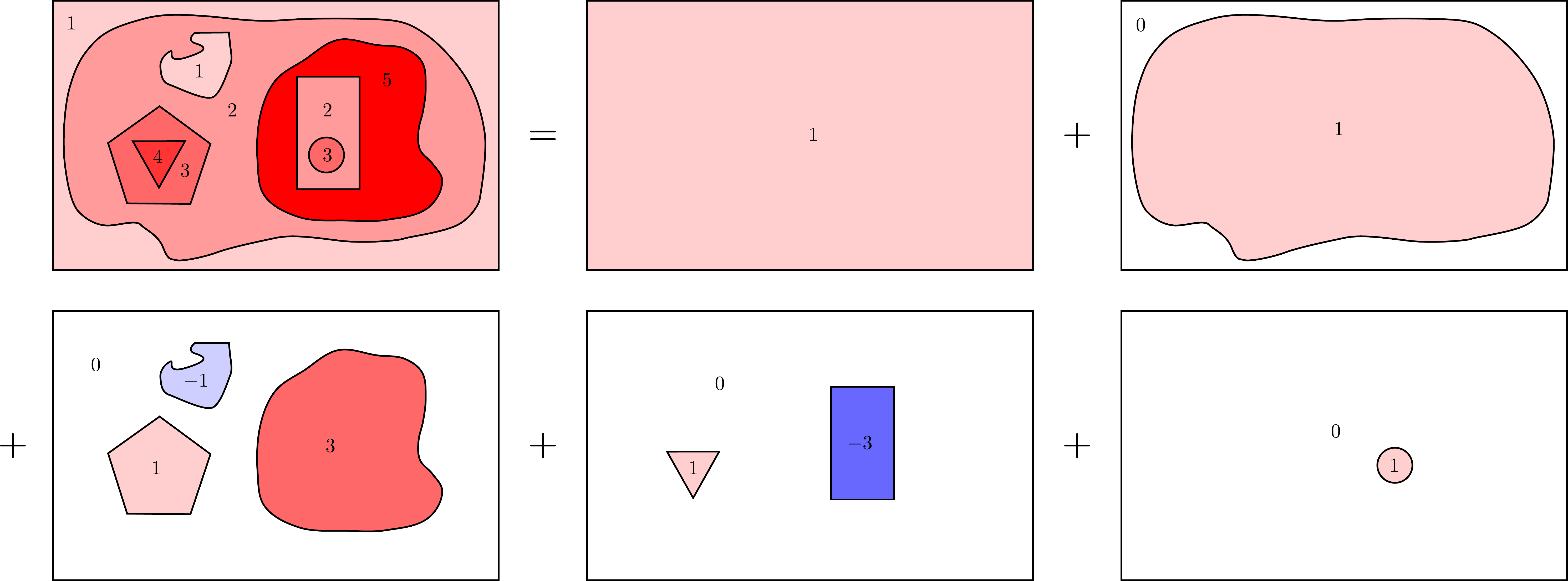}
	\caption{Example of a PCLC conductivity (top left), which is decomposed into a sum of piecewise constant functions on its layers $D_0,\dots,D_4$. The numbers represent function values in each of the colored regions.}
	\label{fig:fig1}
\end{figure}

We need the following notions of a $\tau$-thinning and outer $\tau$-layer of some set $E\subseteq \mathbb{R}^d$ for $\tau>0$:
\begin{align}
	H_\tau(E) &= \{ x\in E \mid \dist(x,\partial E) \geq \tau \}, \label{eq:tauthinning}\\
	F_\tau(E) &= \{ x\in E \mid \dist(x,\partial E) < \tau \}. \label{eq:taulayer}
\end{align}
Next we give some assumptions on sets that will represent material layers for $\gamma$; see also Figure~\ref{fig:fig1}.

\begin{assumption} \label{assump}
	Let $\tau > 0$, $N\in\rum{N}$, and $\{D_j\}_{j=1}^N$ be sets in $\mathbb{R}^d$ satisfying:
	\begin{enumerate}[(i)]
		\item $D_j$ is the closure in $\mathbb{R}^d$ of a nonempty open set with Lipschitz boundary.
		\item $D_j$ has connected complement $\rum{R}^d\setminus D_j$.
		\item $D_{j+1} \subseteq H_\tau(D_j)$ for $j = 1,\dots,N-1$ and $D_1\subset \Omega$.
		\item Each set $D_j$ consists of finitely many connected components $\{D_{j,n}\}_{n=1}^{N_j}$.
	\end{enumerate} 
\end{assumption}

Based on these assumptions, we can define the class of piecewise constant layered conductivities~(PCLC). Let $\chi_E$ denote the characteristic function on $E\subseteq \mathbb{R}^d$.
\begin{definition} \label{def:pclc}
	Suppose $\{D_j\}_{j=1}^N$ satisfy Assumption~\ref{assump} with $\tau>0$, then we call $\gamma$ a PCLC coefficient provided that
	\begin{equation*}
		\gamma = c_0 + \sum_{j=1}^N\sum_{n=1}^{N_j} c_{j,n} \chi_{D_{j,n}},
	\end{equation*}
	where $c_0>0$ and $c_{j,n} \in\rum{R}\setminus\{0\}$ satisfy that $\gamma > 0$ in $\Omega$. Here $D_j$ is called the $j$'th layer of $\gamma$, with $D_0 = \overline{\Omega}$ denoting the $0$'th layer.
	
	For $k\in \{0,1,\dots,N\}$ we define the $k$'th layer-truncated conductivity as
	\begin{equation*}
		\gamma_k = c_0 + \sum_{j=1}^k\sum_{n=1}^{N_j} c_{j,n} \chi_{D_{j,n}}. \label{eq:gammak}
	\end{equation*}
\end{definition}

We now give an algorithm for reconstructing a PCLC coefficient $\gamma$ from its local ND map~$\Lambda(\gamma)$. The method successively reconstructs each layer-truncated conductivity $\gamma_k$, and naturally terminates at $k = N$ for which $\gamma_N = \gamma$. More precisely, the method requires knowledge of the domain and measurement boundary $\Omega$ and $\Gamma$, the datum $\Lambda(\gamma)$, the background conductivity $c_0$ (which may be found via boundary determination, i.e.\ reconstruction of $\gamma$ at the domain boundary), and that $\gamma$ is a PCLC coefficient with some thickness $\tau>0$ between the layers. The method does not require any a priori information about the number of layers $N$, the number of components in each layer $N_j$, or any information on the conductivity values $c_{j,n}$ in each component.

\section{The reconstruction method and its proof}

The reconstruction method starts from the exterior, i.e.\ with $\gamma_0 \equiv c_0$ and with $D_0 = \overline{\Omega}$. It will now be shown how to reconstruct $\gamma_{k+1}$ from $\gamma_k$ and $\Lambda(\gamma)$ for any $k\in\{0,1,\dots,N-1\}$. This defines the method for reconstructing $\gamma = \gamma_N$ one layer at a time.

In the following, for self-adjoint operators $A$ and $B$ on $L^2_\diamond(\Gamma)$, the inequality $A\geq B$ means that 
\begin{equation*}
	\inner{(A-B)f,f}\geq 0
\end{equation*}
for all $f\in L^2_\diamond(\Gamma)$. 

For the purpose of shape reconstruction of $D_{k+1}$, define the family of admissible test inclusions inside $D_k$ as
\begin{align*}
	\mathcal{A}_{k} &= \{ C \subseteq D_k \mid C \text{ is the closure of an open set,}  \\
	&\hphantom{{}= \{C \subset D_k \mid{}}\text{has connected complement,} \\
	&\hphantom{{}= \{C \subset D_k \mid{}}\text{and has Lipschitz boundary } \partial C \}.
\end{align*}
For $C\in\mathcal{A}_k$, let $\Lambda_\mu(\gamma_k;C)$ denote the local ND map for the coefficient given by
\begin{equation*}
	\begin{dcases}
		\mu &\text{in } C, \\
		\gamma_k &\text{elsewhere},
	\end{dcases}
\end{equation*}
where $\mu$ may either take the value $0$ or $\infty$. 
From \cite[Theorem~3.7]{Garde2020c}, it follows immediately that
\begin{equation} \label{eq:recon1}
	D_{k+1} = \cap\{ C\in \mathcal{A}_{k} \mid \Lambda_0(\gamma_k;C) \geq \Lambda(\gamma) \geq \Lambda_\infty(\gamma_k;C) \}.
\end{equation}
\begin{remark}
The inequalities in \eqref{eq:recon1} can be tested numerically by simulating the test ND maps using a peeling-type approach similar to the implementation used in \cite{Garde_2019}. It is also possible to formulate the method for reconstructing the parts of $D_{k+1}$ inside each component of $D_k$ separately, as in \cite{Garde_2019b}, if this turns out to be preferable from a numerical point of view.
\end{remark}

Fix a component $D_{k+1,m_0}$ of $D_{k+1}$. We need to reconstruct the constant $c_{k+1,m_0}$ from $\gamma_k$, $\Lambda(\gamma)$, and $D_{k+1}$. Denote by $\Lambda_\mu(\gamma_k;m_0,t)$ the local ND map for the coefficient
\begin{equation*}
	\begin{dcases}
		\mu &\text{in } D_{k+1}\setminus F_\tau(D_{k+1,m_0}), \\
		\gamma_k + t & \text{in } F_\tau(D_{k+1,m_0}), \\
		\gamma_k &\text{elsewhere}.
	\end{dcases}
\end{equation*}
Here we only consider values of $t\in\mathbb{R}$ such that $\gamma_k + t > 0$ in $D_{k+1,m_0}$. Recall the definition of $F_\tau$ in \eqref{eq:taulayer}. In particular, using extreme inclusions makes it possible to focus locally on the outer $\tau$-layer of $D_{k+1,m_0}$, without having to worry about contributions from the other components or from the $(k+2)$'th layer. An application of \cite[Proof of Theorem~3.7]{Garde2020c} entails the following monotonicity relations on the outer $\tau$-layer,
\begin{align}
	t &\geq c_{k+1,m_0} \qquad \text{if and only if} \qquad \Lambda(\gamma) \geq \Lambda_\infty(\gamma_k;m_0,t), \label{eq:recon2a}\\
	t &\leq c_{k+1,m_0} \qquad \text{if and only if} \qquad \Lambda_0(\gamma_k;m_0,t) \geq \Lambda(\gamma). \label{eq:recon2b}
\end{align}
Since $c_{k+1,m_0} \neq 0$, using either of \eqref{eq:recon2a} or \eqref{eq:recon2b} with $t=0$ determines the sign of $c_{k+1,m_0}$. Once the sign is found, the value of $c_{k+1,m_0}$ can be determined via the one-dimensional optimisation problem:
\begin{equation*}
	c_{k+1,m_0} = \begin{cases}
		\min\{ t > 0 \mid \Lambda(\gamma) \geq \Lambda_\infty(\gamma_k;m_0,t)\} & \text{if } c_{k+1,m_0}>0, \\
		\max\{ t < 0 \mid \Lambda_0(\gamma_k;m_0,t) \geq \Lambda(\gamma)\} & \text{if } c_{k+1,m_0}<0.
	\end{cases}
\end{equation*}
By solving such an optimisation problem for each component of $D_{k+1}$ yields $\gamma_{k+1}$. The method can be repeated until we reach $\gamma_N = \gamma$, at which point subsequent uses of the method will result in an empty set from \eqref{eq:recon1}, thereby indicating that the final layer has been reconstructed.
\begin{remark}
	Note that by a simple modification, we can also allow the inner-most layer to be either perfectly insulating or perfectly conducting. 
	
	Another generalisation allows $\gamma$ to only be piecewise constant on some of the outer-most layers. These layers can still be reconstructed using the method, and for the first non-piecewise constant layer the shape can still be found, although it requires a priori knowledge on the number of layers before the coefficient fails to be piecewise constant. This is e.g.\ relevant for reconstruction of the skull's shape and its conductivity in brain imaging.
\end{remark}

\begin{remark}
	There are still some advantages to using the method in \cite{Garde_2019b} in its original form. Firstly, by not using extreme test inclusions allows a formulation with the Fr\'echet derivative of the forward problem, thereby enabling a faster numerical shape reconstruction of the individual components of a layer. Secondly, the method in \cite{Garde2020c} requires Lipschitz boundaries for the inclusions, thus there cannot be cusps on the interior interfaces as in \cite{Garde_2019b}. Although in practical EIT such cusps are not often present.
\end{remark}

\bibliographystyle{plain}

\begin{thebibliography}{10}
	
	\bibitem{Garde2020c}
	V.~Candiani, J.~Dard\'e, H.~Garde, and N.~Hyv{\"o}nen.
	\newblock Monotonicity-based reconstruction of extreme inclusions in electrical
	impedance tomography.
	\newblock {\em SIAM J. Math. Anal.}, 52(6):6234--6259, 2020.
	
	\bibitem{Esposito2021}
	A.~C. Esposito, L.~Faella, G.~Piscitelli, R.~Prakash, and A.~Tamburrino.
	\newblock Monotonicity {P}rinciple in tomography of nonlinear conducting
	materials.
	\newblock {\em Inverse Problems}, 37(4), 2021.
	\newblock Article ID 045012.
	
	\bibitem{Garde_2019b}
	H.~Garde.
	\newblock Reconstruction of piecewise constant layered conductivities in
	electrical impedance tomography.
	\newblock {\em Comm. PDE}, 45(9):1118--1133, 2020.
	
	\bibitem{GardeHyvonen2021b}
	H.~Garde and N.~Hyv\"onen.
	\newblock Reconstruction of singular and degenerate inclusions in
	{C}alder\'on's problem.
	\newblock 2021.
	\newblock Preprint arXiv:2106.07764 [math.AP].
	
	\bibitem{GardeStaboulis_2016}
	H.~Garde and S.~Staboulis.
	\newblock Convergence and regularization for monotonicity-based shape
	reconstruction in electrical impedance tomography.
	\newblock {\em Numer. Math.}, 135(4):1221--1251, 2017.
	
	\bibitem{Garde_2019}
	H.~Garde and S.~Staboulis.
	\newblock The regularized monotonicity method: detecting irregular indefinite
	inclusions.
	\newblock {\em Inverse Probl. Imag.}, 13(1):93--116, 2019.
	
	\bibitem{Harrach_2019}
	B.~Harrach.
	\newblock Uniqueness and {L}ipschitz stability in electrical impedance
	tomography with finitely many electrodes.
	\newblock {\em Inverse Problems}, 35(2), 2019.
	\newblock Article ID 024005.
	
	\bibitem{Harrach10}
	B.~Harrach and J.~K. Seo.
	\newblock Exact shape-reconstruction by one-step linearization in electrical
	impedance tomography.
	\newblock {\em SIAM J. Math. Anal.}, 42(4):1505--1518, 2010.
	
	\bibitem{Harrach13}
	B.~Harrach and M.~Ullrich.
	\newblock Monotonicity-based shape reconstruction in electrical impedance
	tomography.
	\newblock {\em SIAM J. Math. Anal.}, 45(6):3382--3403, 2013.
	
	\bibitem{Harrach15}
	B.~Harrach and M.~Ullrich.
	\newblock Resolution guarantees in electrical impedance tomography.
	\newblock {\em IEEE T. Med. Imaging}, 34(7):1513--1521, 2015.
	
	\bibitem{Ikehata1998a}
	M.~Ikehata.
	\newblock Size estimation of inclusion.
	\newblock {\em J. Inverse Ill-Posed Probl.}, 6(2):127--140, 1998.
	
	\bibitem{Kang1997b}
	H.~Kang, J.~K. Seo, and D.~Sheen.
	\newblock The inverse conductivity problem with one measurement: stability and
	estimation of size.
	\newblock {\em SIAM J. Math. Anal.}, 28(6):1389--1405, 1997.
	
	\bibitem{Tamburrino2002}
	A.~Tamburrino and G.~Rubinacci.
	\newblock A new non-iterative inversion method for electrical resistance
	tomography.
	\newblock {\em Inverse Problems}, 18(6):1809--1829, 2002.
	
\end{thebibliography}

\end{document}